\documentclass[a4paper,10pt]{article}
\usepackage{amsmath,amsthm,amssymb}
\usepackage{graphicx,subfigure}
\usepackage{url}

\usepackage{algorithm2e}
\usepackage{listings}
\usepackage{color}

\usepackage{multicol}
\def\twoplot[#1]#2#3#4#5{
\begin{figure}[h]
\begin{multicols}{2}
\begin{center}
    \includegraphics*[#1]{#2}
    \caption{\label{#2} #4}
\end{center}
\begin{center}
    \includegraphics*[#1]{#3}
    \caption{\label{#3} #5}
\end{center}
\end{multicols}
\end{figure}
}

\usepackage[colorlinks,bookmarksopen,bookmarksnumbered,citecolor=red,urlcolor=red]{hyperref}
\hypersetup{pdftitle={Use of 3D classified topographic data with FullSWOF for high resolution
simulation of a river flood event over a dense urban area},
bookmarks=true,
pdftoolbar=true,
pdfmenubar=true,
pdfauthor={M. Abily, O. Delestre, L. Amosse, N. Bertrand, C. Laguerre, C.-M. Duluc, P.
Gourbesville},
pdfsubject={flood risks, high resolution modelling, DSM, FullSWOF, photogrammetry, Var river,
2D shallow water equations},
pdfcreator={Delestre},
pdfproducer={Abily},
pdfkeywords={flood hazard}{Shallow-Water equation} {Saint-Venant} {parallelization}
{high resolution modelling} {DSM} {surface model} {photogrammetry}
 {FullSWOF} {Var river} {finite volume} {well-balanced} {hydrostatic reconstruction}}

\title{Use of 3D classified topographic data with FullSWOF for high resolution
simulation of a river flood event over a dense urban area}

\author{{M. Abily}\footnote{Polytech Nice Sophia \& I-CiTy, University of Nice Sophia Antipolis,
France, e-mail : abily@polytech.unice.fr}, {O. Delestre}\footnote{Lab. J.A. Dieudonn\'e \& EPU Nice Sophia, University of
 Nice, France, e-mail : delestre@math.unice.fr}, L. Amoss\'e\footnote{Polytech Nice Sophia \& I-CiTy, University of Nice Sophia
 Antipolis, France}, {N. Bertrand}\footnote{Institut for Radioprotection and Nuclear Safety (IRSN), France},
 {C. Laguerre}\footnote{MAPMO, f\'ed\'eration Denis Poisson, University of Orl\'eans, France, e-mail :
 christian.laguerre@math.cnrs.fr},\\ 
 {C.-M. Duluc}\footnote{Institut for Radioprotection and Nuclear Safety (IRSN), France, e-mail : claire-marie.duluc@irsn.fr}
\; and P. Gourbesville\footnote{Polytech Nice Sophia \& I-CiTy, University of Nice Sophia Antipolis,
France}}

\begin{document}
\maketitle

\begin{abstract} High resolution (infra-metric) topographic data, including photogrammetric born $3$D classified data, are becoming
commonly available at large range of spatial extend, such as municipality or industrial site scale. This category of dataset is
promising for high resolution (HR) Digital Surface Model (DSM) generation, allowing inclusion of fine above-ground structures
which might influence overland flow hydrodynamic in urban environment. Nonetheless several categories of technical and numerical
challenges arise from this type of data use with standard $2$D Shallow Water Equations (SWE) based numerical codes.
\newline
FullSWOF (Full Shallow Water equations for Overland Flow) is a code based on $2$D SWE under conservative form. This code relies on
a well-balanced finite volume method over a regular grid using numerical method based on hydrostatic reconstruction scheme.
When compared to existing industrial codes used for urban flooding simulations, numerical approach implemented in FullSWOF allows
to handle properly flow regime changes, preservation of water depth positivity at wet/dry cells transitions and steady state
preservation. FullSWOF has already been tested on analytical solution library (SWASHES) and has been used to simulate runoff and
dam-breaks. FullSWOF’s above mentioned properties are of good interest for urban overland flow.
\newline
Objectives of this study are ($i$) to assess the feasibility and added values of using HR $3$D classified topographic data to model
river overland flow and ($ii$) to take advantage of FullSWOF code properties for overland flow simulation in urban environment.
\newline
A large aerial $3$D classified topographic data gathering campaign has been conducted by Nice Municipality (France) in $2010$.
Accuracy of this classified data is $0.2$ m in both planimetry and altimetry. This data set is available for the low Var river
valley where an extreme flood event occurred in $1994$. $3$D classified data is used to generate different categories of DSMs
and FullSWOF code has been adapted to be used in a river flood condition context.
\newline
Results demonstrate the feasibility and the clear added value of HR topographical data use. Methodology and limits for such an
approach for engineering perspectives are raised up. The study highlights the need of using HR topographical data set to model
flood event in urban areas and FullSWOF performances for such a purpose are highlighted. 
\end{abstract}

\section{Introduction}\label{sec:intro}
The severity and frequency of urban flooding can be reduced by better planning policies 
 \cite{Djordjevic11}. Numerical modeling tools are commonly used as supporting tools for decision makers willing to assess
 flood mitigation process \cite{Gourbesville09}. At coarse scale, overview of extend and behavior of
 flood event can be estimated by different categories of numerical approaches. These approaches can be based on simplified
 $2$D shallow water equations (SWE) such as diffusive wave or multiple porosity shallow water approaches
 \cite{Guinot12}. Indeed, the cost of these types of numerical approaches, in terms of
computational time and requirement of topographical information, is interesting compared to methods relying on fully resolved $2$D
SWE models. Nevertheless, urbanized areas create complex environment for overland flow locally introducing changes in flow
properties and impacting flow behavior. Therefore, if one objective is to provide a detailed relative comparison of flood hazard
at suburbs and building scale, the use of a fine environment description becomes necessary
 \cite{Sampson12}. For accurate maximal water depth and maximal flow velocity estimations, fully
 resolved $2$D SWE based mathematical model, taking into consideration inertial effects is a more reasonable approach.
 Moreover, the choice of a numerical approach implemented with modern hydrostatic reconstruction and well balanced properties
 will reinforce robustness and accuracy of the computation. FullSWOF\_2D, which relies on finite volume and regular mesh has
 such types of mathematical and numerical properties \cite{Delestre14}. 
\newline
Aerial gathered High Resolution (HR) topographic data is becoming commonly available through specific flight campaign or unmanned
aerial vehicle use (see \cite{Remondino11}). LiDAR and Photogrametry are the most frequently employed
 technologies for this purpose. They allow producing highly accurate DSM, finely describing complex urban environments.
 LiDAR born HR DSM are already commonly used for HR hydraulic modeling. Photogrametry allow through photo-interpretation
 procedure to get accurate and more specifically discriminated (classified) data which allow producing HR DSM adaptable to a
 given hydraulic modeling purpose \cite{Abily14}. The Novelty of using classified data and the
 integration of large amount of data within $2$D hydraulic models are challenging and remain at an experimental stage in
 terms of methodology to optimize HR $3$D classified data use. Feasibility of High Resolution $3$D (HR $3$D) classified
 topographic data use, for hydraulic modeling purpose has been tested at industrial scale for flood event locally generated
 by intense rainfall events \cite{Abily13,Abily14}. The use of such type of data is revealed
 to be both promising and challenging for hydraulic modeling communities.  
\newline
Objectives of this study are ($i$) to assess workability of HR 3D classified topographic data use for river overland flow
modeling and ($ii$) to take advantage of FullSWOF\_2D code properties for overland flow simulations in urban environment. 
\newline
Both FullSWOF\_2D and HR $3$D classified data use for flood river event modeling are tested in this study. A specific approach
has been elaborated, for a medium scale HR DSM creation ($5000$ m per $3500$ m), based on the use of HR $3$D classified data.
Selected area of interest is the low part of the Var river valley. This area has faced, in November $1994$ a flood event. This
area has been covered by a high accuracy photogrammetric data gathering campaign conducted by Nice Municipality (DIGNCA).
 Overland flow influencing structures such as concrete walls, road gutter, sidewalks, etc. are photo-interpreted. These
 structures are included in the dataset and their elevation properties will be a part of the topographic information included
 in the HR DSM specifically created for the hydraulic model. Modifications to FullSWOF\_2D code have been effectuated to model
 river flow scenarios.

\section{Material and methods}\label{sec:mat-meth}

\subsection{High Resolution 3D classified data of the low Var river valley}\label{subsec:High-resol-3D-class}

The photogrammetric data set has been gathered by Nice Municipality Geographic Information Services (DIGNCA) in $2010$-$2011$.
Combination of ($i$) a low altitude flight, ($ii$) a pixel resolution of $0.1$ m at the ground level, ($iii$) a high level
of overlapping among aerial pictures ($80$\%) and ($iv$) the use of an important number of markers for georeferencing
 (about $200$), lead to a high level of accuracy over the urban area of the city. 
\newline
Photo-interpretation allows creation of vectorial information based on photogrammetric dataset 
 \cite{Egels04,Linder06}. A photo-interpreted dataset is composed of classes of points, polylines and polygons
digitalized based on photogrammetric data. Important aspects in the photo-interpretation process are ($i$) classes definition
 and ($ii$) techniques and dataset quality used for photo-interpretation. Both will impact the design of the output
 classified dataset \cite{Lu07}.
Class definition step has to be elaborated prior to the photo-interpretation step. The number, the nature and criteria for
classes’ definition will depend on the objectives of the photo-interpretation campaign. Photo-interpretation techniques can
be made ($i$) automatically by algorithm use, ($ii$) manually by a human operator on a Digital Photogrammetric Workstation
 (DPW) or ($iii$) by a combination of the two methods. The level of accuracy is higher when the photo-interpretation is done by
 a human operator on a DPW, but is much more resources and time consuming \cite{Lafarge10}. 
\newline
Principle of Nice city $3$D HR classified dataset created from photogrammetry is explained in \cite{Andres12}.
 The mean accuracy of the photo-interpreted data over the low Var valley area is $0.15$ - $0.2$ m in both vertical and
 horizontal dimensions. Error in photo-interpretation is estimated to be around $5$\%. These levels of errors and accuracy
 have been checked through terrestrial topographic measurements effectuated by DIGNCA over $10$\% of the domain covered by the
 photogrammetric campaign. The number of class of elements created as vectorialized features is about $50$. The high level
 of accuracy has allowed to photo-interpret thin above ground features as narrow as concrete walls and road gutters. Over
 the part of the low Var river area selected for the study, total number of polyline features represents more than
 $1\,100\,000$ objects introduced under vector form. 

\subsection{FullSWOF}\label{subsec:Fullswof}

FullSWOF stands for Full Shallow Water equations for Overland Flow
(for more details see \cite{Delestre14}). It is a set of open source C++ (ANSI) codes, freely available to the community from the website
 \url{https://sourcesup.renater.fr/projects/fullswof-2d/}. It is distributed under a GPL like free software license. The
 structure of the code is made to facilitate the development of new evolutions. This software resolves the shallow water
 equations thanks to a well-balanced finite volume method based on the hydrostatic reconstruction (introduced in
 \cite{Audusse04c}). This numerical method has good properties: water mass conservation, well-balancedness
 (at least preservation of lake at rest equilibrium) and positivity water height preservation.
The shallow water system in $2$D (SW$2$D) writes:
\begin{equation}
 \left\{\begin{array}{l}
\partial_t{h}+\partial_x(hu)+\partial_y(hv)=0,\\
\partial_t(hu)+\partial_y(hu^2+gh^2/2)+\partial_y(huv)=gh({S_0}_x-{S_f}_x),\\
\partial_t(hv)+\partial_y(huv)+\partial_y(hv^2+gh^2/2)=gh({S_0}_y-{S_f}_y),
        \end{array}
\right.
\end{equation}
where the unknowns are the velocities $u(x,y,t))$ and $v(x,y,t)$ $\left[\text{m}/\text{s}\right]$ and the water height
$h(x,y,t)$ $\left[m\right]$. The subscript $x$ (respectively $y$) stands for the $x$-direction (resp. the $y$-direction):
$-{S_0}_x=\partial_x z(x,y)$ and $-{S_0}_y=\partial_y z(x,y)$ are the ground slopes and ${S_f}_x$ and ${S_f}_y$ the friction terms.
In FullSWOF, we have chosen to solve the SW$2$D on a structured grid. So we have chosen a numerical method adapted to the shallow
water system in $1$D (SW$1$D) and then it is generalized to $2$D thanks to the method of lines. So in what follows, we describe
the numerical method for the SW$1$D. The SW$1$D writes:
\begin{equation}
\left\{\begin{array}{l}
 \partial_t h+\partial_x (hu)=0,\\
 \partial_t (hu)+\partial_x(hu^2+gh^2/2)=gh(S_0-S_f). 
       \end{array}\right.
\end{equation}
 in what follows, we consider Manning's friction law 
 \begin{equation}
S_f=n^2\dfrac{u|u|}{h^{4/3}}=n^2\dfrac{q|q|}{h^{10/3}},  
 \end{equation}
 with $q=hu$ the discharge in $\left[\text{m}^2/\text{s}\right]$. The hydrostatic reconstruction is based on a general principle
 of reconstruction. We begin with a first order finite volume scheme for the form of SW$1$D (without source terms): choosing
 a positive and consistent numerical flux ${\bf F}(U_L,U_R)$ ({\it e.g.} Rusanov, HLL, kinetic, ...), a finite volume scheme
 writes under the general form
\begin{equation}
 \dfrac{U_i^*-U_i^n}{\Delta t}+\dfrac{{\bf F}(U_{i},U_{i+1})-{\bf F}(U_{i-1},U_{i})}{\Delta x}=0,
\end{equation}
where $\Delta t$ is the time step and $\Delta x$ the space step. The idea is to modify this scheme by applying the flux to
reconstructed variables. Reconstruction can be used to get higher order schemes (MUSCL, ENO, ...), in that case higher order in
time is obtained through TVD-Runge-Kutta methods. And the aim of the hydrostatic reconstruction is to be well-balanced. It is
designed to preserve at least steady states at rest ($u=0$). When it is directly applied on the initial scheme, it leads to a
order one scheme, while coupling it with high order reconstruction increases the order and the accuracy.
\newline
We describe now the implementation of this method for high order accuracy. The first step consists in performing a high order
reconstruction (MUSCL, ENO, ...). To treat properly the topography source term $\partial_x z$, this reconstruction is applied on
$u$, $h$ and $h+z$. This gives us the reconstructed variables $(U_-,z_-)$ and $(U_+,z_+)$, on which the hydrostatic reconstruction
is applied 
\begin{equation}
 \left\{\begin{array}{l}
h_{i+1/2L}=\max(h_{i+1/2-}+z_{i+1/2-}-\max(z_{i+1/2-},z_{i+1/2+}),0),\\
U_{i+1/2L}=(h_{i+1/2L},h_{i+1/2L}u_{i+1/2-}),\\
h_{i+1/2R}=\max(h_{i+1/2+}+z_{i+1/2+}-\max(z_{i+1/2-},z_{i+1/2+}),0),\\
U_{i+1/2R}=(h_{i+1/2R},h_{i+1/2R}u_{i+1/2+}).
        \end{array}
\right.
\end{equation}
The finite volume scheme is modified as follows
\begin{equation}
\dfrac{U_i^*-U_i^n}{\Delta t}+\dfrac{F_{i+1/2L}^n-F_{i-1/2R}^n-Fc_i^n}{\Delta x}=0,
\end{equation}
where
\begin{equation}
 F_{i+1/2L}^n=F_{i+1/2}^n+S_{i+1/2L}^n,\quad  F_{i-1/2R}^n=F_{i-1/2}^n+S_{i-1/2R}^n
\end{equation}
are left (resp. right) modifications of the numerical flux for the homogeneous system. In this formula, the flux is now applied
 with reconstructed variables $F_{i+1/2}^n={\bf F}(U_{i+1/2L}^n,U_{i+1/2R}^n)$ and we take
\begin{equation}
\left\{\begin{array}{c}
S_{i+1/2L}^n=\left(\begin{array}{c}
                    0\\
		    \dfrac{g}{2}(h_{i+1/2-}^2-h_{i+1/2L}^2)
                   \end{array}\right),\\
                   \\
S_{i-1/2R}^n=\left(\begin{array}{c}
                    0\\
		    \dfrac{g}{2}(h_{i-1/2+}^2-h_{i-1/2R}^2)
                   \end{array}\right). 
\end{array}\right.
\end{equation}
Finally, for consistency and well-balancing, a centered source term is added
\begin{equation}
 Fc_i=\left(\begin{array}{c}
             0\\
	    -g\dfrac{h_{i-1/2+}+h_{i+1/2-}}{2}(z_{i+1/2-}-z_{i-1/2+})
            \end{array}\right).
\end{equation}
The chosen numerical strategy consists in the HLL flux (see \cite{Delestre10b}, not detailed here) and
 a modified MUSCL reconstruction. It has shown to be the best compromise between accuracy, stability and CPU time cost
 (in \cite{Delestre10b}).
 The MUSCL reconstruction of a scalar variable $s\in \mathbb{R}$ writes
\begin{equation}
 s_{i-1/2+}=s_i-\Delta x .{Ds_i}/2,\quad s_{i+1/2-}=s_i+\Delta x .{Ds_i}/2,
\end{equation}
with the minmod slope limiter
\begin{equation}
Ds_i=minmod\left(\dfrac{s_i-s_{i-1}}{\Delta x},\dfrac{s_{i+1}-s_i}{\Delta x}\right),
\end{equation}
with
\begin{equation}
minmod(x,y)=\left\{\begin{array}{ll}
                           \min(x,y) & \text{if}\,x,y\geq 0,\\
			  \max(x,y) & \text{if}\,x,y\leq 0,\\
			  0 & \text{else}.
                          \end{array}\right. 
\end{equation}
In order to keep the discharge conservation, the reconstruction of the velocity has to be modified as
\begin{equation}
 u_{i-1/2+}=u_i-\dfrac{h_{i+1/2-}}{h_i}\dfrac{\Delta x}{2}Du_i
\quad
u_{i+1/2-}=u_i+\dfrac{h_{i-1/2+}}{h_i}\dfrac{\Delta x}{2}Du_i
\end{equation}

If we take $Ds_i=0$, we recover the first order scheme in space. The friction term is taken into account by a fractional step,
with the following system
\begin{equation}
 \partial_t U=\left(\begin{array}{c}
                     0\\
		    -ghS_f
                    \end{array}\right).
\end{equation}
This system is solved thanks to a semi-implicit method (as in \cite{Delestre10b})
\begin{equation}
 \left\{\begin{array}{l}
         h^{n+1}=h^*,\\
	\dfrac{q^{n+1}-q^*}{\Delta t}=-n^2\dfrac{q^{n+1}|q^n|}{h^n (h^{n+1})^{4/3}}.
        \end{array}\right.
\end{equation}
This method allows to preserve stability (under a classical CFL condition) and steady states at rest. Finally a TVD-Runge Kutta
method is applied to get second order in time. For the generalization to $2$D, we use the HLLC flux introduced in combined with
the method of lines. Concerning boundary conditions, we have modified the code, in order to have the discharge only in the
riverbed, it is based on Riemann invariants. Finally, as we aim at simulating with big data, we have used a parallel version of
FullSWOF based on a domain decomposition and the MPI library \cite{Cordier13}.

\subsection{Method for HR 3D classified data use for hydraulic modeling}\label{Method-HR-3D}

\begin{figure}[htbp]
\begin{center}
\includegraphics[width=0.92\textwidth]{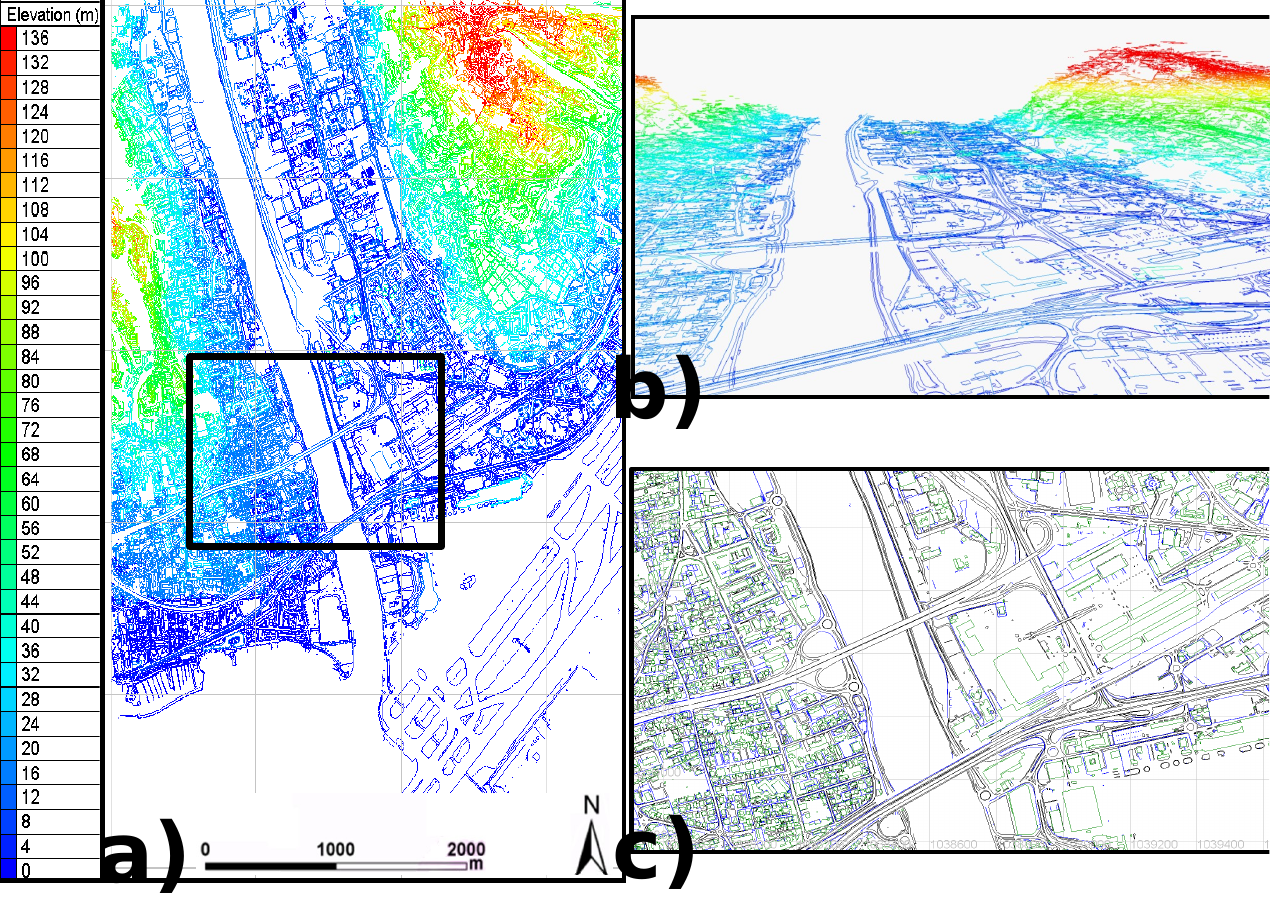}
\caption{Overview (a) and zoom (b) and (c) of the HR $3$D dataset selected classes at step two of the HR DSM creation before
 bridges and flow blocking macro-structures removal.}
\label{fig1}
\end{center}
\end{figure}

To create the HR DSM, the following approach has been carried out. First, a DTM using multiple ground level information
sources: points, polygons and polylines is created and provided at a $0.5$ m resolution by DIGNCA. Then, a selection procedure
among classified data is performed. This selection is achieved by considering concrete elements which can influence overland
flow drainage path only. It includes dikes, buildings, walls and “concrete” above ground elements (such as sidewalks, road
gutters, roundabound, doors steps, etc.). $12$ classes are selected among the $50$ classes of the $3$D photo-interpreted dataset
(figure \ref{fig1}). During this step, polylines giving information on elevated roads and bridges, which might block overland
flow paths, are removed. The remaining total number of polylines is $52\,600$ after these two selection steps. 
\newline

\begin{figure}[htbp]
\begin{center}
\includegraphics[width=0.98\textwidth]{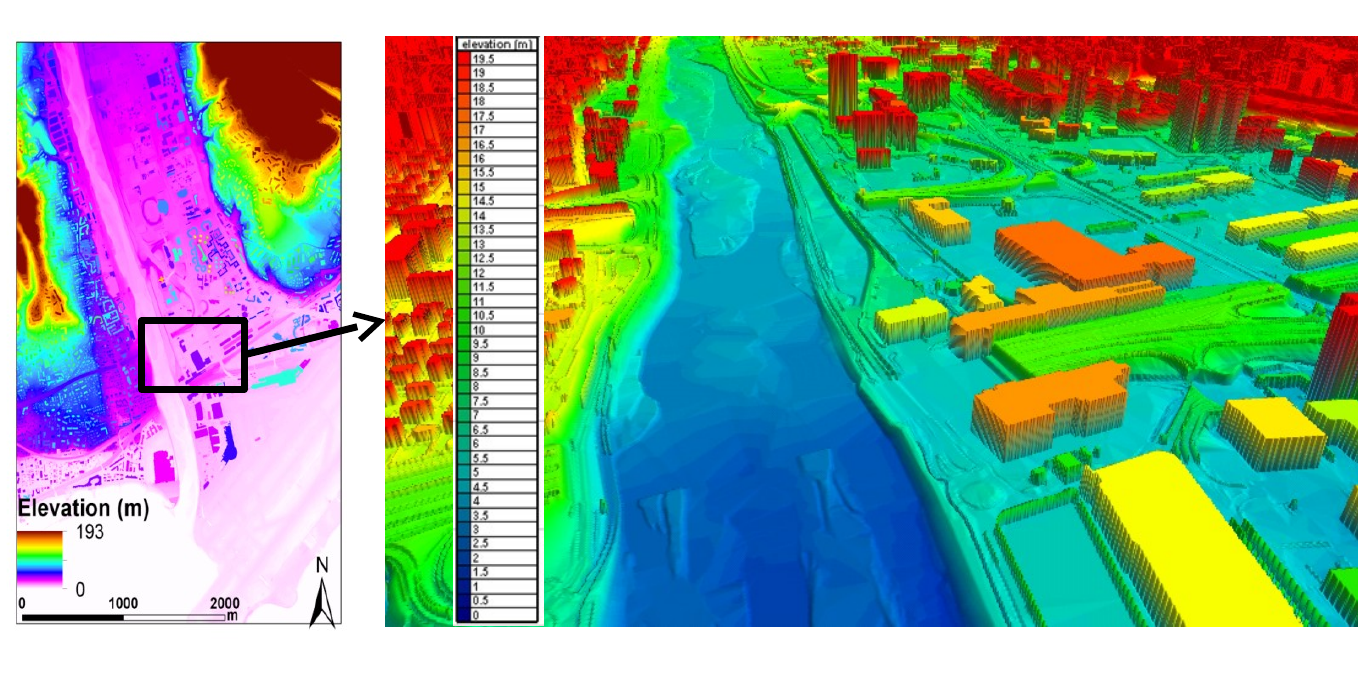}
\caption{HR DSM overview illustrating ground above ground elements elevation (with $z$ axis scale multiplied by $2$ for 
 clarity of the $3$D rendering.}
\label{fig2}
\end{center}
\end{figure}

Final step of HR DSM elaboration consists in extruding elevation information of selected polylines on the DTM. To proceed,
features represented by closed lines are converted to polygons ({\it e.g.} buildings, round abound, sidewalks). Polylignes
 and polygons are then converted to raster at desired resolution (here $1$ m resolution) for extrusion over the DTM. Eventually,
 HR DSM which has elevation information of selected $3$D classified features is produced (figure \ref{fig2}). The HR DSM resolution
 is here $1$ m. This choice of resolution is explained as follow:  it will allow to integrate directly the HR DSM at desired
 regular mesh resolution in FullSWOF\_2D. At this resolution the number of mesh cells is above $17.8$ millions. The previously
 described method has allowed inclusion of thin elements impacting flow behavior of infra-metric dimension, oversized to
 metric size, in the $1$ m resolution regular mesh.

\subsection{Site and river flood event scenario}\label{subsec:site-river-flood}

The $5^\text{th}$ to the $6^\text{th}$ of November $1994$, an important flood event has occurred in the low part of the Var
catchment. This historical flood event had severe consequences. The flood scenario for our tests is based on estimated hydrogram
of this event
\cite{Guinot03}. Through our tests, we want to produce a HR map of maximal water
 depths reached in the low Var valley, using the produced HR DSM with FullSWOF code. Objective here is not to reproduce the flood
 event. Indeed, the site changed a lot since $1994$: levees, dikes and urban structures have been intensively constructed in this
area. To shorten the simulation length, we chose to simulate a $9$ hours scenario. First, a constant discharge of
$1500\,\text{m}^3.\text{s}^{-1}$ is run for $3$ hours to reach a steady state. Then the overtopping part of the hydrogram
is run, reaching the estimated pic discharge ($3700\,\text{m}^3.\text{s}^{-1}$) and then decreasing long enough to observe
a diminution of the overland flow water depth. The Manning-Strickler $n$ coefficient is spatially uniform in overland flow
areas. No energy loss properties have been included in the hydraulic model to represent the bridges piers effects.

\section{Results and discussion}

\begin{figure}[htbp]
\begin{center}
\includegraphics[width=0.98\textwidth]{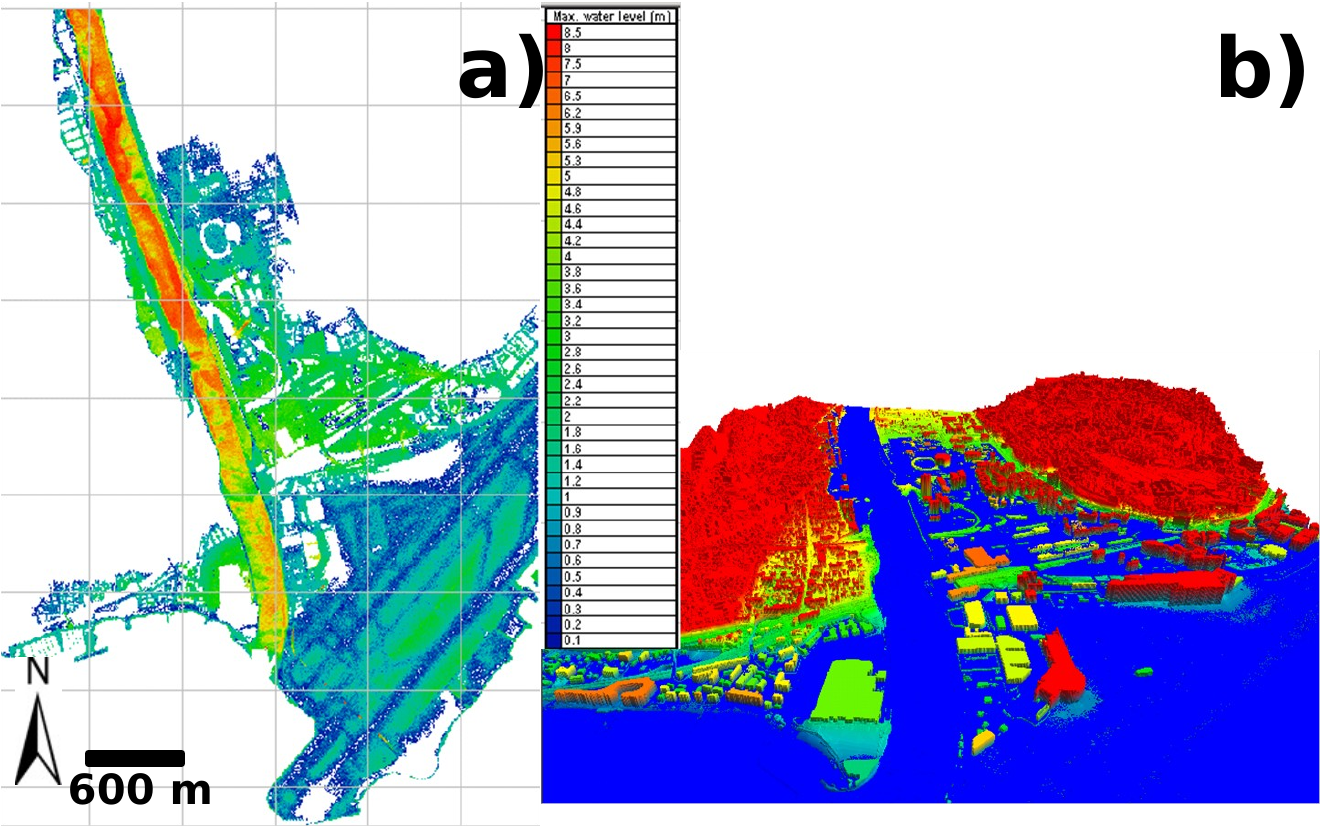}
\caption{Overland flow maximal water depths flood map calculated using HR DSM with FullSWOF\_2D (a) and $3$D global
 representation of flood extent (b).}
\label{fig3}
\end{center}
\end{figure}

An overview of produced overland flow calculated maximal water depths is given in figure \ref{fig3}. High resolution modeling
of flood river event in an urban environment is of a great interest for urban planners as it allows producing detailed maps
of maximal water elevations and maximal velocities maps for a given flood scenario. A proof of concept of $3$D HR classified
data use for river flood modeling is given here. Advantages of such an approach rely: ($i$) possibility to include detailed
surface elements influencing overland flow, and ($ii$) in automatization and modularity of class selection for HR DSM production.
\newline

Two limits, in our HR DSM created for our test study, have to be emphasized. ($i$) The riverbed section itself was filled by
$0.1$ m to $0.2$ m of water at the time of the photogrammetric campaign. Therefore the sections of the river are here
 underestimated, not without standing uncertainties the fact that changes in riverbed occur during a flood event. ($ii$)
 Bridges piers, reducing river section are not included in the HR DSM for our workability test.
\newline
More generally speaking, several categories recommendation and limits deserved to be put to the light for practical
engineering applications using such type of data. Even though HR DSM has a high level of accuracy, HR $3$D data have different
types of inherent errors. It includes white noise, biased, and punctual errors in measurement. Second type of errors is related
to photo-interpretation. It can be omission, addition or misclassification of elements. These types of errors can lead to
important changes in overland flow path in the hydraulic calculation. Moreover, criteria for photo-interpretation have to be
checked as well as they can highly influence HR DSM. For instance a criteria (for human operators or algorithm) might be to
close a polyligne, if distance between two points is below a certain threshold. This will block entrance of water in a given
area. Moreover classification criteria for a given photo-interpreted data set might not have been created specifically for
water modeling purpose. For instance what will be classified, as concrete wall, not be based on material criteria but on
structure width/elevation ratio. In that case permeable structures, such as fences, can be classified as walls. Finally, a
limitation appears regarding bridges piers where information is not given by aerial techniques.
\newline
The $3$D HR classified data are heavy and their manipulation for pre- and post-process is computational resources demanding.
HR DSM use, with fully resolved $2$D SWE codes at this scale, requires use of intensive calculation resources.

\section{Conclusions}

A proof of concept of High Resolution (HR) $3$D classified data use, to produce a HR DSM for river flood simulation in complex
environment has been presented in this study. Hydraulic modeling has been performed by adapting and using FullSWOF which is
a code relying on fully resolved $2$D Shallow Water Equations. Interest has been keen on FullSWOF as his numerical properties of
mass conservation, well-balancedness and positivity preservation are relevant for HR overland flow modeling in urban areas.
\newline
A Method to design a HR DSM including elements influencing overland flow has been presented. For such a purpose, workability of
HR 3D classified topographic data use is relevant at this scale sub-city scale. Existing limits in this approach are put to the
light. These limits mainly consists in ($i$) difficulty to handle this important amount of data, ($ii$) existence of unavoidable
errors of classification in photo-interpretation, ($iii$) classification procedure which might not have been specifically
designed criteria for hydraulic purpose. 

\section*{Acknowledgments}

Photogrametric and photo-interpreted dataset used for this study have been kindly provided by DIGNCA for research purpose.
Technical expertise on DIGNCA dataset has been provided by G. Tacet and F. Largeron. Computations have been performed at
the M\'esocentre d'Aix-Marseille Universit\'e as well as IDRIS resources. Technical support for codes adaptation on high
performance computation centers has been provided by F. Lebas, H. Coullon and P. Navarro.


\begin{thebibliography}{}

\bibitem[Abily et~al., 2014]{Abily14}
Abily, M., Duluc, C.-M., and Gourbesville, P. (2014).
\newblock Use of standard 2d numerical modeling tools to simulate surface
  runoff over an industrial site: Feasibility and comparative performance
  survey over a test case.
\newblock In Gourbesville, P., Cunge, J., and Caignaert, G., editors, {\em
  Advances in Hydroinformatics}, Springer Hydrogeology, pages 19--33. Springer
  Singapore.

\bibitem[Abily et~al., 2013]{Abily13}
Abily, M., Gourbesville, Andres, L., and Duluc, C.-M. (2013).
\newblock Photogrammetric and {L}i{D}{A}{R} data for high resolution runoff
  modeling over industrial and urban sites.
\newblock In Zhaoyin, W., Lee, J. H.-w., Jizhang, G., and Shuyou, C., editors,
  {\em Proceedings of the 35th {I}{A}{H}{R} World Congress, September 8-13,
  2013, Chengdu, China}. Tsinghua University Press, Beijing.

\bibitem[Andres, 2012]{Andres12}
Andres, L. (2012).
\newblock L'apport de la donn\'ee topographique pour la mod\'elisation 3d fine
  et classifi\'ee d'un territoire.
\newblock {\em Revue {X}{Y}{Z}}, 133 - 4e trimestre:24--30.

\bibitem[Audusse et~al., 2004]{Audusse04c}
Audusse, E., Bouchut, F., Bristeau, M.-O., Klein, R., and Perthame, B. (2004).
\newblock A fast and stable well-balanced scheme with hydrostatic
  reconstruction for shallow water flows.
\newblock {\em SIAM J. Sci. Comput.}, 25(6):2050--2065.

\bibitem[{Cordier, S.} et~al., 2013]{Cordier13}
{Cordier, S.}, {Coullon, H.}, {Delestre, O.}, {Laguerre, C.}, {Le, M. H.},
  {Pierre, D.}, and {Sadaka, G.} (2013).
\newblock Fullswof paral: Comparison of two parallelization strategies (MPI and
  SkelGIS) on a software designed for hydrology applications.
\newblock {\em ESAIM: Proc.}, 43:59--79.

\bibitem[Delestre, 2010]{Delestre10b}
Delestre, O. (2010).
\newblock {\em Simulation du ruissellement d'eau de pluie sur des surfaces
  agricoles/ rain water overland flow on agricultural fields simulation}.
\newblock PhD thesis, Universit\'e d'Orl\'eans (in French), available from TEL:
  tel.archives-ouvertes.fr/INSMI/tel-00531377/fr.

\bibitem[Delestre et~al., 2014]{Delestre14}
Delestre, O., Cordier, S., Darboux, F., Du, M., James, F., Laguerre, C., Lucas,
  C., and Planchon, O. (2014).
\newblock Fullswof: A software for overland flow simulation.
\newblock In Gourbesville, P., Cunge, J., and Caignaert, G., editors, {\em
  Advances in Hydroinformatics}, Springer Hydrogeology, pages 221--231.
  Springer Singapore.

\bibitem[Djordjevi\'c et~al., 2011]{Djordjevic11}
Djordjevi\'c, S., Butler, D., Gourbesville, P., Mark, O., and Pasche, E.
  (2011).
\newblock New policies to deal with climate change and other drivers impacting
  on resilience to flooding in urban areas: the {CORFU} approach.
\newblock {\em Environmental Science \& Policy}, 14(7):864--873.
\newblock Adapting to Climate Change: Reducing Water-related Risks in Europe.

\bibitem[Egels and Kasser, 2004]{Egels04}
Egels, Y. and Kasser, M. (2004).
\newblock {\em Digital Photogrammetry}.
\newblock Taylor \& Francis.

\bibitem[Gourbesville, 2009]{Gourbesville09}
Gourbesville, P. (2009).
\newblock Data and hydroinformatics: new possibilities and challenges.
\newblock {\em Journal of Hydroinformatics}, 11(3-4):330--343.

\bibitem[Guinot, 2012]{Guinot12}
Guinot, V. (2012).
\newblock Multiple porosity shallow water models for macroscopic modelling of
  urban floods.
\newblock {\em Advances in Water Resources}, 37(0):40--72.

\bibitem[Guinot and Gourbesville, 2003]{Guinot03}
Guinot, V. and Gourbesville, P. (2003).
\newblock Calibration of physically based models: back to basics?
\newblock {\em Journal of Hydroinformatics}, 5(4):233--244.

\bibitem[Lafarge et~al., 2010]{Lafarge10}
Lafarge, F., Descombes, X., Zerubia, J., and Pierrot~Deseilligny, M. (2010).
\newblock {Structural approach for building reconstruction from a single DSM}.
\newblock {\em Trans. on Pattern Analysis and Machine Intelligence},
  32(1):135--147.

\bibitem[Linder, 2006]{Linder06}
Linder, W. (2006).
\newblock {\em Digital Photogrammetry: A Practical Course}.
\newblock Springer Verlag.

\bibitem[Lu and Weng, 2007]{Lu07}
Lu, D. and Weng, Q. (2007).
\newblock A survey of image classification methods and techniques for improving
  classification performance.
\newblock {\em International Journal of Remote Sensing}, 28(5):823--870.

\bibitem[Remondino et~al., 2011]{Remondino11}
Remondino, F., Barazzetti, L., Nex, F., Scaioni, M., and Sarazzi, D. (2011).
\newblock U{A}{V} photogrammetry for mapping and 3{D} modeling -- {C}urrent
  status and future perspectives.
\newblock In {\em Archives of Photogrammetry, Remote Sensing and Spatial
  Information Sciences}, volume 38(1/C22). ISPRS Conference UAV-g, Zurich,
  Switzerland.

\bibitem[Sampson et~al., 2012]{Sampson12}
Sampson, C.~C., Fewtrell, T.~J., Duncan, A., Shaad, K., Horritt, M.~S., and
  Bates, P.~D. (2012).
\newblock Use of terrestrial laser scanning data to drive decimetric resolution
  urban inundation models.
\newblock {\em Advances in Water Resources}, 41(0):1--17.

\end{thebibliography}

\end{document}